\documentclass[11pt]{amsart}

\usepackage[margin=1in]{geometry}

\usepackage[T1]{fontenc}
\usepackage[latin1]{inputenc}

\usepackage{amsthm, amsmath, amssymb,booktabs,xcolor,graphicx}
\usepackage{mathtools}
\usepackage[shortlabels]{enumitem}
\usepackage{scalerel,stackengine}
\usepackage[textsize=scriptsize,colorinlistoftodos]{todonotes}
\usepackage[hidelinks]{hyperref}
\usepackage[noabbrev,capitalize,nameinlink]{cleveref}
\crefformat{equation}{(#2#1#3)}
\crefrangeformat{equation}{(#3#1#4) to~(#5#2#6)}
\crefmultiformat{equation}{(#2#1#3)}
{ and~(#2#1#3)}{, (#2#1#3)}{ and~(#2#1#3)}

\usepackage[color,final]{showkeys}

\colorlet{refkey}{orange!50}
\colorlet{labelkey}{blue!50}

\newtheorem{theorem}{Theorem}[section]

\newtheorem{prop}[theorem]{Proposition}
\newtheorem{lemma}[theorem]{Lemma}

\newtheorem{conj}[theorem]{Conjecture}
\newtheorem*{question*}{Question}

\theoremstyle{definition}

\newtheorem*{definition*}{Definition}

\theoremstyle{remark}
\newtheorem{rem}[theorem]{Remark}
\newtheorem{alg}[theorem]{Algorithm}

\numberwithin{equation}{section}

\newcommand{\abs}[1]{\left\lvert#1\right\rvert}

\newcommand{\floor}[1]{\left\lfloor #1 \right\rfloor}
\newcommand{\ceil}[1]{\left\lceil #1 \right\rceil}
\newcommand{\paren}[1]{\left( #1 \right)}
\newcommand{\sqb}[1]{\left[ #1 \right]}
\newcommand{\set}[1]{\left\{ #1 \right\}}

\newcommand{\eps}{\varepsilon}

\renewcommand{\P}{\mathbf{P}}

\stackMath
\newcommand\widecheck[1]{
\savestack{\tmpbox}{\stretchto{
  \scaleto{
    \scalerel*[\widthof{\ensuremath{#1}}]{\kern-.6pt\bigwedge\kern-.6pt}
    {\rule[-\textheight/2]{1ex}{\textheight}}
  }{\textheight}
}{0.5ex}}
\stackon[1pt]{#1}{\scalebox{-1}{\tmpbox}}
}
\DeclareMathOperator{\Bin}{Bin}

\newcommand{\RR}{\mathbb{R}}
\newcommand{\PP}{\mathbb{P}}

\newcommand{\mcC}{\mathcal{C}}
\newcommand{\cF}{\mathcal{F}}
\newcommand{\mcF}{\mathcal{F}}

\newcommand{\mcI}{\mathcal{I}}

\newcommand{\mcN}{\mathcal{N}}

\makeatletter
\newcommand\thankssymb[1]{\textsuperscript{\@fnsymbol{#1}}}
\makeatother

\title{On the number of error correcting codes}

\author{Dingding Dong}
\author{Nitya Mani}
\author{Yufei Zhao}

\thanks{Mani was supported by the NSF Graduate Research Fellowship Program and a Hertz Graduate Fellowship.}
\thanks{Zhao was supported by NSF CAREER award DMS-2044606, a Sloan Research Fellowship, and the MIT Solomon Buchsbaum Fund.}

\address{Dong: Department of Mathematics, Harvard University, Cambridge, MA 02138, USA}
\email{ddong@math.harvard.edu}

\address{Mani \& Zhao: Department of Mathematics, Massachusetts Institute of Technology, Cambridge, MA 02139, USA}
\email{\{nmani,yufeiz\}@mit.edu}

\begin{document}

\begin{abstract}
We show that for a fixed $q$, the number of $q$-ary $t$-error correcting codes of length $n$ is at most $2^{(1 + o(1)) H_q(n,t)}$ for all $t \leq (1 - q^{-1})n - C_q\sqrt{n \log n}$ (for sufficiently large constant $C_q$),
where $H_q(n, t) = q^n / V_q(n,t)$ is the Hamming bound and $V_q(n,t)$ is the cardinality of the radius $t$ Hamming ball.
This proves a conjecture of Balogh, Treglown, and Wagner, who showed the result for $t = o(n^{1/3} (\log n)^{-2/3})$.
\end{abstract}

\maketitle

\section{Introduction}

A \textit{$q$-ary $t$-error correcting code of length $n$} is a subset of $[q]^n$ whose elements pairwise differ in at least $2t+1$ coordinates. Error correcting codes play a central role in coding theory, as they allow us to correct errors in at most $t$ coordinates when sending codewords over a noisy channel. In the course of developing the container method, Sapozhenko in~\cite{SAP05} posed the following natural question: how many error correcting codes (with given parameters) are there? 

Given a $t$-error correcting code $\mcC \subseteq [q]^n$ of maximum size, any subset of $\mcC$ is also a $t$-error correcting code. 
Consequently, the number of $q$-ary  $t$-error correcting codes of length $n$ is at least $2^{\abs{\mcC}}$. 
It is natural to suspect that the total number of $t$-error correcting codes is not much larger, provided that this maximum size code $\mcC$ is not too small.

A fundamental problem in coding theory is determining the maximum size of an error correcting code with given parameters (e.g., see \cite[Chapter 5]{vL99} and \cite[Chapter 17]{MS77} for classic treatments). 
There is an easy geometric upper bound on the size of a $q$-ary $t$-error correcting code. 
We use $d(x,y)$ to denote the \emph{Hamming distance} between $x$ and $y \in [q]^n$, defined as the number of coordinates on which $x$ and $y$ differ.
Let $B_q(x, r) := \set{ y \in [q]^n : d(x,y) \le r}$
be the \textit{Hamming ball} of radius $r$ centered at $x \in [q]^n$ with associated cardinality
\[
V_q(n, r) := \sum_{i = 0}^r \binom{n}{i}(q-1)^{i}.
\]
Any $t$-error correcting code in $[q]^n$ has size at most
\[
H_q(n, t) := \frac{q^n}{V_q(n, t)}.
\]
This simple upper bound is known as the 
\textit{Hamming bound} (also known as the \textit{sphere packing bound}).
Codes that attain equality for this upper bound are called perfect codes; the only non-trivial perfect codes are Hamming codes and Golay codes (see \cite{vL75}).  Except for some special parameters such as these, there is a large gap between the best known upper and lower bounds.

The Hamming bound can be significantly improved when $t \ge C \sqrt{n}$ (see~\cite{Del73,MRRW77,vL99}). Most studies on bounds for codes focus on the linear distance case, i.e., $d = \floor{\delta n}$ for some constant $0 < \delta < 1$.
Curiously, even after consulting with several experts on the subject, we are not aware of any upper bound better than the Hamming bound when $t < c\sqrt{n}$. 
An important family of constructions is BCH codes, which can achieve sizes roughly $H_q(n, t)/t!$; this is within a subexponential but superpolynomial factor of the Hamming bound.
As we are primarily interested in the regime $t \le c\sqrt{n}$ in this paper, we view the Hamming bound $H_q(n,t)$ as a proxy for the size of the code.

Balogh, Treglown and Wagner~\cite{BTW16} showed that as long as $t = o(n^{1/3} (\log n)^{-2/3})$, 
the number of binary $t$-error correcting codes of length $n$ is at most $2^{(1+o(1)) H_2 (n,t)}$ (their paper only considered the $q=2$ case, though the situation for larger fixed $q$ is similar).
They conjectured that the same upper bound should also hold for a larger range of $t$.
We prove a stronger version of this conjecture.

\begin{theorem}\label{t:main}
Fix integer $q \ge 2$. Below $t$ is allowed to depend on $n$.
\begin{enumerate}
    \item [(a)] For $0 < t \le 10\sqrt{n}$, the number of $q$-ary $t$-error correcting  codes of length $n$ is at most $2^{(1+o(1)) H_q (n,t)}$.
    \item [(b)] There exists some constant $C_q > 0$ such that for all $10\sqrt{n} < t \leq (1 - 1/q) n - C_q\sqrt{n \log n}$, the number of $q$-ary $t$-error correcting  codes of length $n$ is $2^{o(H_q(n,t))}$.
    \item [(c)] There exists some constant $c_q > 0$ such that for all $t > (1-1/q) n - c_q \sqrt{n \log n}$, one has $H_q(n,t) = O(n^{1/10})$.
    (Since there are $q^n$ one-element codes, the number of codes is not $2^{O(H_q(n,t))}$ whenever $H_q(n, t) = O(n^{1/10})$.)
\end{enumerate}
\end{theorem}

The most interesting part is (a), where $t \le 10\sqrt{n}$; 
here, our proof builds on the container argument of \cite{BTW16}. 
The key idea in \cite{BTW16} was analyzing the container algorithm in two stages depending on the number of remaining vertices.
Our work introduces two new ideas. We first give a pair of more refined supersaturation estimates; these estimates arise from understanding Hamming ball intersections and weak transitivity-type properties of a Hamming distance graph we study. Secondly, we employ a variant of the classical graph container algorithm with early stopping.
See~\cref{sec:container} for the container argument as well as the statements of the required supersaturation results, and~\cref{sec:supersat} for proofs of the supersaturation estimates.

For the remaining cases when $t > 10 \sqrt{n}$, we observe in~\cref{sec:bounds} that the Elias bound (a classic upper bound on code sizes) implies the desired result when $t \le (1 - q^{-1})n - c_q \sqrt{n \log n}$ and that $H_q(n, t)$ is very small when $t$ is any larger.
Also, note that by the Plotkin bound \cite[(5.2.4)]{vL99}, if $2t+1 - (1-1/q)n = \alpha > 0$, then the size of a $t$-error correcting codes in $[q]^n$ is at most $(2t+1)/\alpha$, which is very small, and so (b) is uninteresting when $t$ is large.

A $q$-ary code of \emph{length} $n$ and \emph{distance} $d$ is a subset of $[q]^n$ whose elements are pairwise separated by Hamming distance at least $d$. The maximum size of such a code is denoted $A_q(n,d)$.
It is natural to conjecture the following much stronger bound on the number of codes of length $n$ and distance $d$, which is likely quite difficult.

\begin{conj}\label{c:stronger}
For all fixed $c > 0$, the number of $q$-ary codes of length $n$ and distance $d$ is $2^{O(A_q(n,d))}$ whenever $d < (1 - q^{-1} - c)n$.
\end{conj}

The main difficulty seems to be a lack of strong supersaturation bounds, e.g., showing that sets with size much larger than $A_q(n,d)$ must have many pairs of elements with distance at most $d$.

The above conjecture is analogous to the following conjecture of Kleitman and Winston~\cite{KW82}.
We use $\operatorname{ex}(n,H)$ to denote the maximum number of edges in an $n$-vertex $H$-free graph.

\begin{conj}
	For every bipartite graph $H$ that contains a cycle, the number of $H$-free graphs on $n$ labeled vertices is $2^{O_{H}(\operatorname{ex}(n, H))}$.
\end{conj}

Kleitman and Winston~\cite{KW82} proved the conjecture for $H = C_4$, and in doing so developed what is now called the graph container method, which is central to the rest of this paper. 
See Ferber, McKinley, and Samotij~\cite{FMS20} for recent results and discussion on the problem of enumerating $H$-free graphs.

Let $r_k(n)$ denote the maximum size of a $k$-AP-free subset of $[n]$. Balogh, Liu, and Sharifzadeh~\cite{BLS17} recently proved the following theorem. 

\begin{theorem}[Balogh, Liu, and Sharifzadeh~\cite{BLS17}]
Fix $k \ge 3$. The number of $k$-AP-free subsets of $[n]$ is $2^{O(r_k(n))}$ for infinitely many $n$.	
\end{theorem}

A notable feature of the above theorem is that the asymptotic order of $r_k(n)$ is not known for any $k \ge 3$, similar to the situation for $A(n,d)$.  
Ferber, McKinley, and Samotij~\cite{FMS20} proved analogous results for the number of $H$-free graphs.
It remains an open problem to extend the above result to all $n$.
Also, Cameron and Erd\H{o}s asked whether the number of $k$-AP-free subsets of $[n]$ is $2^{(1+o(1))r_k(n)}$. This is unknown, although for several similar questions the answer is no (see~\cite{BLS17} for discussion).
Likewise, one can ask whether the number of $q$-ary error-correcting codes of length $n$ and distance $d$ is $2^{(1+o(1))A_q(n,d)}$.

\textit{Remarks about notation.} For nonnegative quantities depending on $n$, we write $f \lesssim g$ to mean that there is some constant $C$ such that $f \le Cg$ for all sufficiently large $n$. Throughout the paper we view $q$ as a constant, and so hidden constants are allowed to depend on $q$.

Given a finite set $X$, we write $2^X$ for the collection of all subsets of $X$, 
$\binom{X}{t}$ for the collection of $t$-element subsets of $X$,
and $\binom{X}{\le t}$ for the collection of subsets of $X$ with at most $t$ elements.
We also write $\binom{n}{\le t} = \sum_{0 \le i \le t} \binom{n}{i}$.

\section{Graph container argument} \label{sec:container}

In this section, we analyze the  case $t \le 10\sqrt{n}$.
As in~\cite{BTW16}, we use the method of graph containers.
This technique was originally introduced by Kleitman and Winston~\cite{KW80,KW82}. 
See Samotij~\cite{Sam15} for a modern introduction to the graph container method.
We will not need the more recent hypergraph container method for this work.
Instead, as in \cite{BTW16}, we refine the graph container algorithm and analysis.

Let $G_{q,n,t}$ be the graph with vertex set $[q]^n$ 
where two vertices are adjacent if their Hamming distance is at most $2t$.
For a subset $S \subseteq [q]^n,$ we denote by $G_{q,n,t}[S]$ the induced subgraph of $G_{q, n,t}$ on $S$.
For simplicity of notation we write $G[S] = G_{q,n,t}[S]$, as the tuple $(q,n,t)$ does not change in the proof.
Given a graph $G$, 
we write $\Delta(G)$ for the maximum degree of $G$
and $i(G)$ for the number of independent sets in $G$.

Applications of the container method usually require supersaturation estimates. Our main advance, building on \cite{BTW16}, 
is a more refined supersaturation bound for codes. 
The first lemma below already allows us to prove the main theorem for $t = o(\sqrt{n}/\log n)$; in concert with the second lemma, we are able to establish the full theorem.

\begin{lemma}[Supersaturation I]\label{lem:n32t12}
Let $t \le 10 \sqrt{n}$. 
If $S \subseteq [q]^n$ has size $|S| \ge n^4H_q(n,t)$, then
\[
\Delta(G[S]) \gtrsim \frac{n^{3/2}}{ H_q(n,t)}|S|.
\]
\end{lemma}

\begin{lemma}[Supersaturation II] \label{lem:degis} 
Let $\eps>0$. For all sufficiently large $n$, if $60 \le t\le 10\sqrt{n}$ and $S \subseteq[q]^n$ satisfies $\Delta(G[S]) \le n^5$, then $i(G[S]) \le 2^{(1+\eps)H_q(n,t)}$.
\end{lemma}

We defer the proof of the above two lemmas to \cref{sec:sstrong}. 
Assuming these two lemmas, we now run the container argument to deduce the main theorem.
The container algorithm we use is a modification of the classical graph container algorithm of Kleitman and Winston~\cite{KW80,KW82}. As in \cite{BTW16},  we divide the container algorithm into several stages and analyze them separately. Moreover, we allow ``early stopping'' of the algorithm if we can certify that the residual graph has few independent sets.

\begin{lemma}[Container] \label{lem:container}
For every $\eps > 0$, the following holds for $n$ sufficiently large. For all $60 \le t \le 10 \sqrt{n}$, there exists a collection $\mcF$ of subsets of $[q]^n$ with the following properties:
\begin{itemize}
    \item $|\mcF| \le 2^{\eps H_q(n,t)}$;
    \item Every $t$-error correcting code in $[q]^n$ is contained in $S$ for some $S \in \mcF$;
    \item $i(G[S]) \le 2^{(1 + \eps )H_q(n, t)}$ for every $S \in \mcF$.
\end{itemize}
\end{lemma}

This container lemma implies~\cref{t:main}(a).

\begin{proof}[Proof of~\cref{t:main}(a)]\label{p:thma}
If $t < 60$, then the result is already known (\cite{BTW16} proved it for $t = o(n^{1/3} (\log n)^{-2/3})$).
So assume $60 \le t \le 10\sqrt{n}$.
Let $\eps > 0$. 
Thus, letting $\mcF$ be as in Lemma~\ref{lem:container}, the number of $t$-error correcting codes in $[q]^n$ is at most 
\[
\sum_{S \in \mcF} i(G[S]) \le |\mcF| \cdot 2^{(1 + \eps) H_q(n,t)} = 2^{(1 + 2\eps) H_q(n, t)},
\]
for sufficiently large $n$. 
Since $\eps > 0$ can be arbitrarily small, we have the desired result.
\end{proof}

\begin{proof}[Proof of \cref{lem:container}]
Fix $\eps > 0$.  Let $C > 0$ be a sufficiently large constant, and assume that $n$ is sufficiently large throughout.
We show that there exists some function 
\[
f: \binom{V(G)}{\le \frac{C H_q(n, t) \log n}{n}} \rightarrow 2^{V(G)},
\]
such that for every independent set $I \in \mcI(G)$, there is some ``fingerprint'' $P \subseteq I$ that satisfies the following three conditions:
\begin{enumerate}
    \item $\abs{P} \le \frac{C H_q(n, t) \log n}{n}$,
    \item $i(G[P\cup f(P)])\le 2^{(1+\eps)H_q(n,t)}$, 
    \item $I \subseteq P \cup f(P)$.
\end{enumerate}
We consider the following algorithm that constructs $P$ and $f(P)$ given some $I \in \mcI(G).$

\begin{alg}
Fix an arbitrary order $v_1, \ldots, v_{q^n}$ of the elements of $V(G)$. Let $I \in \mcI(G)$ be an independent set, as the input.
\begin{enumerate}
    \item Initialize $G_0 := G$ and $P:= \emptyset$.
    \item For $i = 1, 2, \ldots$
    \begin{itemize}
        \item If $i(G_{i-1})\le 2^{(1+\eps)H_q(n,t)}$, then set $f(P)=V(G_{i-1})$ and terminate.
        \item Else, let $u$ be a maximum degree vertex in $G_{i-1}$ (break ties by taking the earliest vertex in the given ordering).
        \item If $u \not \in I$, define $G_i$ by taking $G_{i-1}$ and removing $u$, and continue onto the next step.
        \item If $u \in I$, then add $u$ to $P$, define $G_i$ by taking $G_{i-1}$ and removing $u$ and its neighbors, and continue.
    \end{itemize}
     \item Output $P$ and $f(P)$.
\end{enumerate}
\end{alg}
Given any independent set $I$, the algorithm outputs some $P$ and $f(P)$. Note that $f(P)$ depends only on $P$ (i.e., two different $I$'s that output the same $P$ will always output the same $f(P)$). Also, for any independent set $I$, the output satisfies $P \subseteq I \subseteq P \cup f(P)$. 

To prove $\abs{P} \le \frac{C H_q(n, t) \log n}{n}$, we distinguish stages in the algorithm based on the size of $V(G_i)$.
\begin{itemize}
    \item Let $P_1$ denote the vertices added to $P$ when $|V(G_{i-1})| \ge n^4 H_q(n, t)$.
    \item Let $P_2$ denote the vertices added to $P$ when $|V(G_{i-1})| < n^4 H_q(n, t)$.
\end{itemize}
While we are adding vertices to $P_1$, according to \cref{lem:n32t12}, we are removing at least a $\beta \gtrsim n^{3/2}/H_q(n,t)$ fraction of vertices with every successful addition. We have
\begin{align*}
    |P_1|
    \le  \frac{\log\left( \frac{q^n}{n^4H_q(n,t)}\right)}{\log\left(\frac{1}{1-\beta}\right)}
    &\lesssim \frac{\log V_q(n,t)-4\log n}{\beta}
     \le \frac{\log ( (nq)^t)}{\beta}
    \lesssim \frac{t\log n }{n^{3/2}} H_q(n,t)
    \lesssim \frac{ \log n}{n} H_q(n, t).
\end{align*}
While we are adding vertices to $P_2$, by~\cref{lem:degis}, we are removing at least $n^5$ vertices with every successful addition, as $i(G_{i-1}) \ge 2^{(1+\eps)H_q(n,t)}$ during the associated iteration. Consequently, 
\begin{align*}
    |P_2|\lesssim \frac{n^4 H_q(n,t)}{n^5}= \frac{H_q(n,t)}{n}.
\end{align*}
These two inequalities imply that
\[
|P| = \abs{P_1} + \abs{P_2} \le  \frac{C \log n}{n} H_q(n, t).
\]
Thus
\[
    i(G[P\cup f(P)])
    \le 2^{\abs{P}}\cdot  i(f(P))
    \le 2^{\eps H_q(n, t) } \cdot 2^{(1+\eps)H_q(n,t)}\le 2^{(1+2\eps)H_q(n,t)}.
\]

Take $\cF$ to be the collection of all $P\cup f(P)$ obtained from this procedure as $I$ ranges over all independent sets of $G$. We have just proved that  every $t$-error correcting code is contained in some $S\in\cF$. We also have $i(G[S])\leq 2^{(1+2\eps)H_q(n,t)}$ for every $S\in\cF$ (equivalent to the stated result after replacing $\eps$ by $\eps/2$). Finally, using that $\log \binom{n}{m} \le m \log(O(n)/m)$, we have
\begin{align*}
    \log \abs{\cF}
    &\le \log \binom{q^n}{\le \frac{C \log n}{n} H_q(n, t) }
    \\
    &\lesssim \frac{C \log n}{n}  H_q(n, t)  \log \left(O\paren{\frac{n V_q(n,t)}{C\log n}}\right)
    \\
    &\le \frac{C \log n}{n} \log\left( O\paren{\frac{n (qn)^t}{C\log n}}\right) H_q(n, t) 
    \\
    &\lesssim \frac{C t (\log n)^2}{n} H_q(n,t).
\end{align*}
Since $t \le 10 \sqrt{n}$, $\abs{\cF} \le 2^{\varepsilon H_q(n,t)}$ for sufficiently large $n$.
\end{proof}

It remains to prove the supersaturation estimates \cref{lem:n32t12,lem:degis}, which we will do in~\cref{sec:supersat} after first observing some technical preliminaries in~\cref{sec:ballestimates}. The cases when $t > 10\sqrt{n}$ are in~\cref{sec:bounds}.

\section{Hamming ball volume estimates}\label{sec:ballestimates}

\subsection{Hamming ball volume ratios}
We first record some estimates about sizes of Hamming balls of different radii, beginning with the following decay estimate.

\begin{lemma}\label{lem:vntdecreasing}
For every $1 \le i\le t$, we have 
\[
V_q(n-i,t-i)\le \left(\frac{t}{(q-1)n}\right)^i V_q(n,t).
\]
\end{lemma}
\begin{proof}
Since $\binom{n-1}{j}=\frac{j+1}{n} \binom{n}{j+1}$, we have
\begin{align*}
    \frac{V_q(n-1,t-1)}{V_q(n,t)}&=\frac{\sum_{j=0}^{t-1}\binom{n-1}{j}(q-1)^j}{1+\sum_{j=0}^{t-1}\binom{n}{j+1}(q-1)^{j+1}}\le \frac{\sum_{j=0}^{t-1}\binom{n-1}{j}(q-1)^{j}}{\sum_{i=0}^{t-1}\binom{n}{j+1}(q-1)^{j+1}}\le \frac{t}{(q-1)n}.
\end{align*}
Consequently,
\begin{align*}
    \frac{V_q(n-i,t-i)}{V_q(n,t)}&=\frac{V_q(n-1,t-1)}{V_q(n,t)}\cdot\dots \cdot\frac{V_q(n-i,t-i)}{V_q(n-i+1,t-i+1)}\\
    &\le \frac{1}{(q-1)^i}\cdot\frac{t}{n}\cdot\frac{t-1}{n-1}\cdot\dots\cdot \frac{t-i+1}{n-i+1}\le \left(\frac{t}{(q-1)n}\right)^i. 
    \qedhere
\end{align*}
\end{proof}

\begin{lemma}\label{l:alphavbound}
For $1 \le \alpha \le t$, we have
\[
V_q(n, t + \alpha) \ge \left( \frac{(q-1)n}{t+\alpha}\right)^{\alpha} \left(\frac{n - \alpha + 1 - t}{n - \alpha + 1}\right)^{\alpha} V_q(n, t).
\]
\end{lemma}
\begin{proof}
Observe that 
\[
\frac{V_q(n-\alpha, t)}{V_q(n, t)} = \frac{\sum_{i = 0}^t (q-1)^i {n- \alpha \choose i}}{\sum_{i = 0}^t (q-1)^i {n \choose i}} \ge \left(\frac{n - \alpha + 1 - t}{n - \alpha + 1}\right)^{\alpha}.
\]
Also, by~\cref{lem:vntdecreasing}, we have
\[
V_q(n-\alpha, t) \le \left(\frac{t+\alpha}{(q-1)n}\right)^{\alpha} V_q(n, t+ \alpha).
\]
Combining this pair of inequalities gives
\[
V_q(n, t+ \alpha) \ge \left( \frac{(q-1)n}{t+\alpha}\right)^{\alpha} V_q(n-\alpha, t) \ge \left( \frac{(q-1)n}{t+\alpha}\right)^{\alpha} \left(\frac{n - \alpha + 1 - t}{n - \alpha + 1}\right)^{\alpha} V_q(n, t). \qedhere
\] 
\end{proof}

In the regime of large $t$, we will be able to observe~\cref{t:main}(c), i.e., that there exists $c_q > 0$ such that for all $t > (1-1/q) n - c_q \sqrt{n \log n}$, we have $H_q(n,t) = O(n^{1/10})$.
\begin{proof}[Proof of~\cref{t:main}(c)]
 Let $t=(1-1/q)n-\alpha$ with $\alpha<c_q\sqrt{n\log n}$ for some sufficiently small $c_q > 0$. 
 Let $\theta=1-q^{-1}$.
 Let $Y\sim \Bin(n,\theta)$ be a binomial random variable and $Z\sim \mcN(0,1)$ a standard normal random variable. 
 Observe that for all $0\leq i\leq n$, we have $V_q(n, t) = q^n \P(Y \le t)$. 
 By the Barry--Esseen theorem (a quantitative central limit theorem), for all $x\in\RR$,
 \begin{align*}
    \abs{\PP\sqb{\frac{Y-\theta n}{\sqrt{\theta(1-\theta)n}}\leq x}-\PP\sqb{Z\leq x}}\lesssim \frac{1}{\sqrt{n}}.
\end{align*}
We have the following standard estimate on the Gaussian tail:
\[
\P[Z \le  -x] > \frac{1}{\sqrt{2\pi}} \frac{x}{x^2+1} e^{-x^2/2} \qquad \text{for all }x>0.
\]
Thus
\begin{align*}
    \frac{1}{H_q(n,t)} 
    = 
    \frac{V_q(n,\theta n-\alpha)}{q^n} 
    &\geq \P\sqb{Z\leq \frac{-c_q\sqrt{n\log n}}{\sqrt{\theta (1-\theta )n}}} -O\left( \frac{1}{\sqrt{n}}\right) \\
   &\gtrsim  e^{- \Theta(c_q^2 \log n) }-O\left( \frac{1}{\sqrt{n}}\right) \\
    &\gtrsim  n^{-1/10} 
\end{align*}
by choosing $c_q$ to be sufficiently small.
Therefore $H_q(n,t) = O(n^{1/10})$.
\end{proof}

\begin{rem}
It is possible to give more precise estimates of binomial tails with larger deviations.
The above application of the Barry--Esseen theorem is a concise way to obtain what we need.
\end{rem}

\subsection{Intersection of Hamming balls}\label{sec:wtsk}
We record some estimates about the sizes of intersections of Hamming balls as a function of the distance between their centers. Throughout this subsection, we assume $t \le 10\sqrt{n}.$

Let $W_q (n,t, k)$ be the size of the intersection of two Hamming balls in $[q]^n$ of radius $t$, the centers being distance $k$ apart.
It is easy to check that
\begin{align*}
    W_q(n,t,k)&=\sum_{r=0}^k\sum_{s=0}^{k-r}\binom{k}{r} \binom{k-r}{s}(q-2)^{k-r-s}V_q(n-k,t-\max\{k-r,k-s\}).
\end{align*}
Indeed, suppose the two balls are centered at $(1^n)$ and $(2^k1^{n-k})$.
The intersection then consists of all points with $r$ $1$'s and $s$ $2$'s among the first $k$ coordinates (as $r$ and $s$ range over all possible values),
and $\le t - \max\{k-r,k-s\}$ non-$1$ coordinates among the remaining $n-k$ coordinates.

In particular, it is easy to see that 
\[
W_q(n,t,1) = q V_q(n-1,t-1).
\]
(As above, it does not matter what happens in the first coordinate and the remaining $n-1$ coordinates can have at most $t-1$ non-$1$ coordinates). 

\begin{lemma} \label{lem:wdecreasing}
$W_q(n,t, k+1) \leq W_q(n,t, k)$ for every integer $0 \le k \le t$.
\end{lemma}
\begin{proof}
It suffices to consider radius $t$ Hamming balls centered at 
$u = 1^n$,
$v = 2^{k}1^{n-k}$,
and
$w = 2^{k+1}1^{n-k-1}$.
We  show that $\abs{B_q(u,t)\cap B_q(w,t)} \leq \abs{B_q(u,t)\cap B_q(v,t)}$.
Let us consider the following subsets of $[q]^n$:
$$X = B_q(u, t) \cap B_q(w,t) \setminus B_q(v,t) \quad \text{ and } \quad  Y = B_q(u, t) \cap B_q(v,t) \setminus B_q(w,t).$$
Observe that every $x = (x_1, \ldots, x_n) \in X$ must have $x_{k+1} = 2$ and $d(x, w) = t$, $d(x, v) = t+1$.
Define $\phi:X \rightarrow Y$ where $\phi(x)$ is obtained from $x$ by changing the $(k+1)$st coordinate from $2$ to $1$. This map is well-defined since for every $x \in X$, $d(\phi(x), u) < d(x, u) \le t$, $d(\phi(x), w) = t + 1$ and $d(\phi(x), v) = t$, and so $\phi(x) \in Y$. 
Since $\phi$ is injective, we have $|X|\leq|Y|$ and therefore $\abs{B_q(u,t)\cap B_q(w,t)} \leq \abs{B_q(u,t)\cap B_q(v,t)}$.
\end{proof}

\begin{lemma}\label{lem:wsizebound}
For every integer $k\ge 0$, we have
\[
W_q(n,t, 2k+2) 
\leq W_q(n,t, 2k+1) 
\le 2\left(\frac{q^2t}{(q-1)n}\right)^{k}
W_q(n,t,1)
\le
2\left(\frac{q^2t}{(q-1)n}\right)^{k+1}V_q(n,t).
\]
\end{lemma}
\begin{proof}
The first inequality is simply \cref{lem:wdecreasing}.
We further have that
\begin{align*}
    \MoveEqLeft W_q(n,t, 2k+1)
    \\
    &=\sum_{r=0}^{2k+1}\sum_{s=0}^{2k+1-r}{2k+1\choose r} {2k+1-r\choose s}(q-2)^{2k+1-r-s}V_q(n-2k-1,t-2k-1+\min\{r,s\})\\
    &\le 2\sum_{r = 0}^{k}\sum_{s=0}^{2k+1-r}{2k+1\choose r} {2k+1-r\choose s}(q-2)^{2k+1-r-s}V_q(n-2k-1,t+r-2k-1)\\
    &\le 2\sum_{r = 0}^{k}\sum_{s=0}^{2k+1-r}{2k+1\choose r} {2k+1-r\choose s}(q-2)^{2k+1-r-s}V_q(n+r-2k-1,t+r-2k-1).
\end{align*}
Using \cref{lem:vntdecreasing} and substituting the above inequality, we deduce that 
\begin{align*}
\frac{W_q(n,t, 2k+1)}{W_q(n,t,1)} 
&=\frac{W_q(n,t, 2k+1)}{qV_q(n-1,t-1)}\\
&\overset{(*)}{\le}  \frac{2}{q}\sum_{r = 0}^{k}\sum_{s=0}^{2k+1-r}{2k+1\choose r} {2k+1-r\choose s}(q-2)^{2k+1-r-s} \left(\frac{t-1}{(q-1)(n-1)}\right)^{2k-r}\\
&\leq \left(\frac{t-1}{(q-1)(n-1)}\right)^{k}\cdot\frac{2}{q}\sum_{r = 0}^{k}\sum_{s=0}^{2k+1-r}{2k+1\choose r} {2k+1-r\choose s}(q-2)^{2k+1-r-s}\\
&= \left(\frac{t-1}{(q-1)(n-1)}\right)^{k}\cdot\frac{2}{q}\sum_{r = 0}^{k}{2k+1\choose r} ((q-2) + 1)^{2k+1-r}\\
&\le \left(\frac{t-1}{(q-1)(n-1)}\right)^{k}\cdot\frac{2}{q}\sum_{r = 0}^{2k+1}{2k+1\choose r} (q-1)^{2k+1-r}\\
&= \left(\frac{t-1}{(q-1)(n-1)}\right)^{k}\cdot\frac{2}{q} \cdot  q^{2k+1} \\
&\leq \left(\frac{t-1}{(q-1)(n-1)}\right)^{k}\cdot 2q^{2k}\le 2\left(\frac{q^2t}{(q-1)n}\right)^k,
\end{align*}
where $(*)$ follows from the previous expansion of $W_q(n,t,2k+1)$ combined with~\cref{lem:vntdecreasing}.
This proves the second inequality in the lemma. 
To prove the last inequality in the lemma, apply \cref{lem:vntdecreasing} again to yield
\[
   \frac{W_q(n,t, 1)}{V_q(n,t)}=\frac{qV_q(n-1,t-1)}{V_q(n,t)}\leq\frac{qt}{(q-1)n}. 
   \qedhere 
\]
\end{proof}

\section{Supersaturation}\label{sec:supersat}

\subsection{Supersaturation I}\label{sec:sstrong}

From now on, let $S \subseteq [q]^n$ and let $G[S] = (S, E)$ be the associated graph with the edge set
\[
E = \{\{x, y\} \in E : d(x, y) \leq 2t\}.
\]
For each $1 \le  k \le 2t$, write
\[
E_k = \{\{x, y\} \in E : d(x, y) = k\}.
\]
In other words, $E$ is the set of pairs in $S$ with Hamming distance at most $2t$, 
and each $E_k$ is the set of pairs in $S$ with Hamming distance exactly $k$.
Also, for each $v \in S$, denote number of edges in $E_k$ incident to $v$ by 
\[
\deg_k (v) = \abs{\set{ u \in S : d(u,v)} = k}.
\]
We have the trivial bound $\deg_k(v) \le {n \choose k}(q-1)^k$ for every $v \in S$, and hence $\abs{E_k} \le \frac{1}{2}\binom{n}{k} (q-1)^k \abs{S}$.

We first recall the following supersaturation estimate from \cite{BTW16} (stated there for binary codes). 
We include the proof below for completeness.
Later on, we will derive new and stronger supersaturation estimates.

\begin{lemma}[{\cite[Lemma 5.3]{BTW16}}]\label{lem:btw16-supersaturation2}
If $|S|\ge 2H_q(n,t)$, then
\[
\sum_{k=1}^{2t}W_q(n,t,k)|E_k|\ge\frac{|S|^2V_q(n,t)^2}{10\cdot q^n}\quad \text{and}\quad |E| \ge \frac{n|S|^2}{20tH_q(n,t)}.
\]
\end{lemma}
\begin{proof}
Define $K_x= \{a \in S \mid d(x, a) \le t\}$ for $x \in [q]^n$. Observe that
\[
\sum_{k = 1}^{2t} W_q(n,t,k) |E_k| = \sum_{x \in [q]^n} {|K_x| \choose 2},
\]
since both terms count pairs $(x, \{a, b\})$ for $x \in [q]^n$ and distinct $a, b \in S$ such that $d(x,a), d(x, b) \le t$.
From \cref{lem:vntdecreasing} and \cref{lem:wdecreasing}, we know that
$$W_q(n, t, k) \le W_q(n,t,1)=qV_q(n-1,t-1) \le \frac{qt}{(q-1)n} V_q(n, t)$$
for $1 \le k \le 2t$. By convexity, since the average value of $K_x$ over $x \in [q]^n$ is $|S|V_q(n, t)/q^n \ge 2$, we have
\[
\frac{qt}{(q-1)n}\cdot V_q(n, t)|E| \ge W_q(n,t, 1)|E| \ge \sum_{k = 1}^{2t} W_q(n,t,k) |E_k| = \sum_{x \in [q]^n} {|K_x| \choose 2} \ge \frac{|S|^2 V_q(n, t)^2}{10 q^n}.
\]
This proves the first inequality stated in the lemma.
Rearranging, we derive the second inequality:
\[
|E| \ge \frac{(q-1)n|S|^2 V_q(n, t) }{10q t\cdot q^n} \ge \frac{n|S|^2 }{20 t H_q(n, t)}. \qedhere 
\]
\end{proof}

Now let us prove the first supersaturation estimate, \cref{lem:n32t12}, which says that if $t \le 10 \sqrt{n}$ and $|S| \ge n^4H_q(n,t)$, then $\Delta(G[S]) \gtrsim \frac{n^{3/2}}{ H_q(n,t)}|S|$.

\begin{proof}[Proof of Lemma~\ref{lem:n32t12}]
By Lemma~\ref{lem:btw16-supersaturation2}, since $|S| \ge n^4H_q(n,t)$, 
\[
|E| \ge\frac{|S|n}{20t H_q(n,t)}\cdot |S|\ge\frac{n^5}{20t}\cdot |S|.
\]
Together with $|E_1|+|E_2| \le n^2 q^2 |S|$ and $|E_3|+|E_4|\le n^4 q^4|S|$, 
we have (recall that hidden constants are allowed to depend on $q$)
\[
|E_1|+|E_2|\lesssim \frac{t}{n^3} |E| 
\quad \text{and} \quad
|E_3|+|E_4|\lesssim \frac{t}{n} |E|.
\]
Furthermore, by \cref{lem:wsizebound}, 
\begin{align*}
W_q(n,t,2) &\le W_q(n,t,1) ,
\\
W_q(n,t,4) &\le W_q(n,t,3) \lesssim \frac{t}{n} W_q(n,t,1) ,   \qquad \text{and}
\\
W_q(n,t,k) &\le W_q(n,t,5) \lesssim \frac{t^2}{n^2} W_q(n,t,1)  \qquad \text{for all $k \ge 5$}.
\end{align*}
Hence we have
\begin{align*}
    \sum_{k=1}^{2t}W_q(n,t,k)|E_k|
    &\le W_q(n,t,1)(|E_1|+|E_2|)+W_q(n,t,3)(|E_3|+|E_4|)+W_q(n,t,5)|E|
    \\
    &\lesssim \paren{\frac{t}{n^3} + \frac{t}{n} \cdot \frac{t}{n} + \frac{t^2}{n^2}} W_q(n,t,1)|E| 
    \lesssim \frac{t^2}{n^2}\cdot W_q(n,t,1)|E|.
\end{align*}
By~\cref{lem:btw16-supersaturation2}, we have
\begin{align*}
    \frac{|S|^2V_q(n,t)^2}{q^n}
    \lesssim   \sum_{k=1}^{2t}W_q(n,t,k)|E_k|
    \lesssim  \frac{t^2}{n^2}\cdot W_q(n,t,1)|E|.
\end{align*}
Rearranging and then applying $W_q(n,t,1) \lesssim (t/n) V_q(n,t)$ from \cref{lem:wsizebound}, and $t \lesssim \sqrt{n}$, we obtain
\begin{align*}
    |E|\gtrsim
    \frac{n^2}{t^2}\cdot \frac{V_q(n,t)}{W_q(n,t,1)} \cdot \frac{V_q(n,t)}{  q^n} \cdot |S|^2
    \gtrsim \frac{n^3}{t^3} \cdot \frac{1}{H_q(n,t)} \cdot |S|^2
    \gtrsim \frac{ n^{3/2}}{H_q(n,t)}\cdot |S|^2.
\end{align*}
Thus the average degree (and hence the maximum degree) in $G[S]$ is $\gtrsim \frac{n^{3/2}}{H_q(n,t)} |S|$.
\end{proof}

\subsection{Supersaturation II}\label{sec:wsupersat}
We maintain the notation from the previous subsection.

\begin{lemma}\label{lem:maxdegree-vs-indepsets-1}
Suppose $t\le 10\sqrt{n}$ and $\Delta(G[S])\le n^5$.
Fix $\eps > 0$. 
Let
\[
S_1 = \set{ v\in S : \deg_k(v) \le \eps n^{\ceil {k/2}/2} \text{ for each }k=1,\dots,20}.
\]
Then $|S_1|\le (1 + O(\eps))H_q(n,t)$.
\end{lemma}

\begin{proof}
The idea is that the Hamming balls $B_q(v,t)$, $v \in S_1$, are ``mostly disjoint'', in the sense that the overlap is negligible.

For any $v\in S_1$, the overlap of $B_q(v,t)$ with other balls $B_q(u,t)$, $u \in S_1$, has size
\begin{align*}
    \abs{B_q(v,t)\cap \bigcup_{u\in S_1\setminus\{v\}} B_q(u,t)}
    &\le \sum_{k=1}^{2t}\deg_k(v) W_q(n,t,k)\\
    &\le \biggl(\sum_{k=1}^{20}\deg_k(v) W_q(n,t,k)\biggr) + n^5 W_q(n,t,21)
\end{align*}
by \cref{lem:wdecreasing}. Recall from \cref{lem:wsizebound} that $W_q(n,t,k) \le 2(q^2t/((q-1)n))^{\ceil{k/2}} V_q(n,t)$, which is $\lesssim n^{-\ceil{k/2}/2} V_q(n,t)$ for each $1 \le k \le 20$. Combining with $\deg_k(v) \le \eps n^{\ceil{k/2}/2}$ for each $v \in S$, the above sum is $\lesssim \eps V_q(n,t)$. 
In other words, for every $v\in S_1$, the ball $B_q(v,t)$  contains $(1-O(\eps)) V_q(n,t)$ unique points not shared by other such balls, and thus the union of these balls has size $\ge (1-O(\eps)) V_q(n,t)\abs{S_1}$. Since the union is contained in $[q]^n$, we deduce $\abs{S_1} \le (1+O(\eps))q^n / V_q(n,t) = (1+O(\eps)) H_q(n,t)$.
\end{proof}

\begin{lemma} \label{lem:maxdegree-vs-indepsets-2}
Suppose $t \ge 60$ and $\Delta(G[S])\le n^{5}$.
Fix $\eps > 0$. Let
\begin{align*}
    S_2 &=
    \set{v\in S: \sum_{k=1}^{20} \deg_k(v) \ge \frac{\log n}{\eps}}.
\end{align*}
Then the number of independent subsets of $S_2$ satisfies
\[
i(G[S_2])\le 2^{O(\eps H_q(n,t))}.
\]
\end{lemma} %\edit{Isn't this actually $o(H_q(n, t))$. Maybe let's make the Theorem~\ref{lem:degis} say that if so?}
\begin{proof} %\edit{Reread and check mathcal font}
Pick a maximum size subset $X \subseteq S_2$ where every pair of distinct elements of $X$ have Hamming distance greater than $t$. 
Using the definition of $S_2$, 
for each $x \in X$, we can find an $\ceil{(\log n)/\eps}$-element subset $A_x \subseteq S \cap B_q(x,20)$.
For distinct $x,y \in X$, one has $d(x,y) > t \ge 60$, and thus $A_x \cap A_y = \emptyset$.

Consider the balls $B_q(u,t)$, $u \in \bigcup_{x \in X} A_x$. 
As in the previous proof, we will show that these balls are mostly disjoint. 
The intersection of one of these balls with the union of all other such balls has size at most 
\[
\frac{\log n}{\eps} W_q(n,t,1) + n^5 W_q(n,t,21) 
\lesssim \paren{\frac{\log n}{\eps}\cdot\frac{t}{n} + n^5 \cdot \frac{t^{11}}{n^{11}}} V_q(n,t)
 = o(V_q(n,t)).
\]
Indeed, for each of the $\ceil{(\log n)/\eps}-1$ points $u'$ that lie in the same $A_x$ as $u$, the overlap is $\le W_q(n,t,1)$;
each ball $B_q(u',t)$ with $u'$ in some $A_y$ other than $A_x$
contributes to $\le W_q(n,t,21)$ overlap
since $d(u,u') \ge d(x,y) - 40 \ge t-39 \ge 21$.
There are at most $\Delta(G[S] \le n^5$ other balls that intersect a given ball.
Here we invoke the monotonicity of $W_q(n,t,\cdot)$ (\cref{lem:wdecreasing}), and also \cref{lem:wsizebound}.

Thus the union of the balls $B_q(u,t)$, $u \in \bigcup_{x \in X} A_x$, has size $\abs{X} \ceil{(\log n)/\eps} \cdot (1-o(1))V_q(n,t)$.
Since this quantity is at most $q^n$,
\[
\abs{X} \le \frac{q^n}{\ceil{(\log n)/\eps} \cdot (1-o(1))V_q(n,t)} \lesssim  \frac{\eps H_q(n,t)}{\log n}. 
\]

Note that $\bigcup_{x \in X} B_q(x,t) = [q]^n$ (since if there were some uncovered $x' \in [q]^n$, then one can add $x'$ to $X$). 
Every pair of elements of $B_q(x,t)$ has Hamming distance at most $2t$.
It follows that every independent set $I$ in $G[S_2]$ (so the elements of $I$ are separated by Hamming distance $>2t$) can be formed by choosing at most one element from $B_q(x,t) \cap S_2$ for each $x \in X$. Since $\abs{B_q(x,t) \cap S_2} \le n^5$, we have
\[
i(G[S_2]) \le (n^5 + 1)^{\abs{X}} \le 2^{O(\eps H_q(n,t))}. \qedhere 
\]
\end{proof}

Now we prove the second supersaturation estimate, \cref{lem:degis}, 
which says that if $60 \le t\le 10\sqrt{n}$ and $\Delta(G[S]) \le n^5$, then $i(G[S]) \le 2^{(1+\eps)H_q(n,t)}$ for all sufficiently large $n$.

\begin{proof}[Proof of \cref{lem:degis}]
Define $S_1$ and $S_2$ as in the previous two lemmas.
One has $S = S_1 \cup S_2$ provided that $n$ is sufficiently large. Thus
\[
i(G[S])
\le
i(G[S_1]) i(G[S_2])
\le 2^{\abs{S_1}} i(G[S_2])
\le 2^{(1+O(\eps)) H_q(n,t)}.
\]
This is equivalent to the desired result after changing $\eps$ by a constant factor.
\end{proof}

\section{Bounds on codes with larger distances}\label{sec:bounds}

When $t > 10\sqrt{n}$, we observe that a maximum sized $t$ error correcting code must be much smaller than the Hamming bound; this implies an associated bound on the number of $t$-error correcting codes.
Recall that $A_q(n, 2t+1)$ is the maximum size of a $t$-error correcting code over $[q]^n$.
The number of $q$-ary $t$-error correcting codes of length $n$ is at most
\begin{equation}\label{eq:weak-ub}
\binom{q^n}{\le A_q(n, 2t+1)} \le q^{n A_q(n, 2t+1)}.
\end{equation}
If $A_q(n, 2t+1) = o(H_q(n, t)/n)$, then the number of $t$-error correcting codes of length $n$ is $2^{o(H_q(n,t))}$. 
This is indeed the case when $t > 10\sqrt{n}$ due to the following classic upper bound on $A_q(n,d)$.

\begin{theorem}[Elias bound {\cite[Theorem 5.2.11]{vL99}}] \label{t:elias}
Let $\theta = 1 - q^{-1}$. For every $r \le \theta n$ satisfying 
\[
r^2 - 2\theta nr + \theta  nd > 0,
\] we have the upper bound
\[
A_q(n, d) \le \frac{\theta  nd}{r^2 - 2\theta nr + \theta  nd} \cdot \frac{q^n}{V_q(n, r)}.
\]
\end{theorem}

We have the following consequence of the above Elias bound.
\begin{prop}\label{c:larget}
Let $\theta = 1-q^{-1}.$
Let $C_q > 0$ be a sufficiently large constant. If $10\sqrt{n} < t \leq  \theta n - C_q\sqrt{n \log n}$, we have that 
$$A_q(n, 2t+1) = o\left( \frac{H_q(n,t)}{n}\right).$$
\end{prop} 
\begin{proof}

We apply the Elias bound with $r = t + \alpha$, where $\alpha = 7$ if $10 \sqrt{n} < t < n^{4/5}$ and $\alpha = \sqrt{n\log n}$ if $n^{4/5} \le t \leq\theta n - C_q\sqrt{n \log n}$. When $t > 10\sqrt{n}$,
we satisfy the condition 
\begin{equation}\label{e:eliasdenom}
r^2 - 2 \theta nr + \theta n(2t+1)=r^2-\theta(2\alpha-1)n \gtrsim r^2.
\end{equation}
Indeed, when $10\sqrt{n} < t < n^{4/5},$ we have $\alpha =7$, and thus $\theta(2\alpha - 1)n\leq 2\alpha n<14 n$ whereas $r^2 > (10\sqrt{n})^2 = 100 n$. When $t \ge n^{4/5}$, we have $r^2 \ge n^{8/5}$ whereas $\theta(2\alpha-1)n=O(n^{3/2}\sqrt{\log n})$.
By~\cref{l:alphavbound}, we have
\[
V_q(n, r) = V_q(n, t+ \alpha) \ge \left( \frac{(q-1)n}{t+\alpha}\right)^{\alpha} \left(\frac{n - \alpha + 1 - t}{n - \alpha + 1}\right)^{\alpha} V_q(n, t).
\]
Applying~\cref{t:elias} gives that 
\begin{align*}
A_q(n, 2t+1) &\le \frac{\theta n(2t+1)}{r^2 - 2\theta nr + \theta n(2t+1)} \cdot \frac{q^n}{V_q(n, r)} \\
&\overset{\mathclap{\eqref{e:eliasdenom}}}\lesssim \frac{nt}{(t + \alpha)^2} \cdot \frac{q^n}{V_q(n, t+\alpha)}. \\
&\le \frac{n}{t} \cdot \left( \frac{(t+\alpha)(n-\alpha+1)}{(q-1)n(n-\alpha + 1-t)}\right)^{\alpha} H_q(n, t).
\end{align*}

When $10\sqrt{n} < t < n^{4/5}$ and $\alpha = 7$, the above upper bound  simplifies to
\[
A_q(n, 2t+1) \lesssim  \frac{n}{t} \cdot \left( \frac{(t+\alpha)(n-\alpha+1)}{(q-1)n(n-\alpha + 1-t)}\right)^{\alpha} H_q(n, t) \lesssim \frac{t^6}{n^{6}} H_q(n, t) \le \frac{1}{n^{6/5}} H_q(n, t) = o\left( \frac{H_q(n, t)}{n} \right).
\]
When $n^{4/5} \le t \le \theta n - C_q\sqrt{n \log n}$ and $\alpha = \sqrt{n\log n}$,
we have
\begin{align*}
    A_q(n, 2t+1) &\le \frac{\theta n(2t+1)}{r^2 - 2\theta nr + \theta n(2t+1)} \cdot \frac{q^n}{V_q(n, r)}\\
    &\lesssim  \frac{n}{t}  \left( \frac{(t+\alpha)(n-\alpha+1)}{(q-1)n(n-\alpha + 1-t)}\right)^{\alpha} H_q(n, t) \\
  %  &\lesssim\frac{n}{t} \left( \frac{\theta n-(C_q-1)\sqrt{n\log n}}{(q-1)(n-\theta n+(C_q-1)\sqrt{n\log n})}\right)^{\alpha}\left(\frac{n-\sqrt{n\log n}}{n}\right)^\alpha H_q(n, t)\\
   % &\leq \frac{n}{t} \left(\frac{\theta n-C_q\sqrt{n\log n}}{\theta n+(q-1)C_q\sqrt{n\log n}}\right)^\alpha \left(\frac{n-\sqrt{n\log n}}{n}\right)^\alpha H_q(n, t)\\
    &= \frac{n}{t}\left(\frac{t+\sqrt{n\log n}}{(q-1)(n-\sqrt{n\log n}+1-t)}\right)^{\sqrt{n\log n}} \left(\frac{n-\sqrt{n\log n}+1}{n}\right)^{\sqrt{n\log n}}H_q(n,t)\\
    % &\leq \frac{n}{t}\left(\frac{\theta n - C_q\sqrt{n \log n}+\sqrt{n\log n}}{(q-1)(n-\sqrt{n\log n}+1-\theta n + C_q\sqrt{n \log n})}\right)^{\sqrt{n\log n}} \left(\frac{n-\sqrt{n\log n}+1}{n}\right)^{\sqrt{n\log n}}H(n,t)\\
    &\lesssim \frac{n}{t}\left(\frac{\theta n - (C_q-1)\sqrt{n \log n}}{(q-1)(q^{-1}n+ (C_q-1)\sqrt{n \log n})}\right)^{\sqrt{n\log n}} \left(\frac{n-\sqrt{n\log n}}{n}\right)^{\sqrt{n\log n}}H_q(n,t)\\
   % &=\frac{n}{t}\left(1-\frac{q(C_q-1)\sqrt{n\log n}}{\theta n+(q-1)(C_q-1)\sqrt{n\log n}}\right)^{\sqrt{n\log n}}\left(\frac{n-\sqrt{n\log n}}{n}\right)^{\sqrt{n\log n}}H_q(n, t)\\
    &\lesssim \frac{n}{t}\left(1-\frac{(C_q-1)\sqrt{\log n}}{\theta\sqrt{n}}\right)^{\sqrt{n\log n}}\left(\frac{n-\sqrt{n\log n}}{n}\right)^{\sqrt{n\log n}}H_q(n, t)\\
    &\leq \frac{n}{t}\cdot n^{-(C_q-1)/\theta}\cdot n^{-1} H_q(n, t)\\
    &=t^{-1}\cdot n^{-(C_q-1)/\theta} H_q(n, t)=o(n^{-1} H_q(n, t)) 
\end{align*}
by picking $C_q$ sufficiently large. This shows the desired result.
%Picking some $\alpha$ such that $\frac{3 \log n}{\log(1 + 7\eps) - \log(1 - 7\eps)}\leq \alpha  = o(t),$ 
% we find per above that 
% \[
% A_q(n, 2t+1) \lesssim \frac{n}{t} \cdot \left( \frac{(1+\eps)t}{n-(1+\eps)t}\right)^{\alpha} H_q(n, t)\leq  \frac{n}{t}  \left( \frac{(1+\eps)(1/2-4\eps)}{1-(1+\eps)(1/2-4\eps)}\right)^{\alpha} H_q(n, t)\lesssim n^{-2} H_q(n, t),
% \]
% (note that when $t = \Theta(\sqrt{n}),$ taking $\alpha = 5$ suffices). 
\end{proof}
% \begin{rem}When $t > n/2$ (which can occur for $q > 2$), the result becomes trivial as $A_q(n,2t + 1)=1$. The same proof strategy shows that by taking $C_q = C(q, \beta) > 0$ sufficiently large, $A_q(n, 2t+1) \lesssim  \frac{1}{n^{\beta}}H_q(n, t)$ for any $\beta > 0$ whenever $10\sqrt{n} < t \le \theta n - C_q \sqrt{n \log n}$.
% \end{rem}

\begin{proof}[Proof of Theorem~\ref{t:main}(b)]
When $10\sqrt{n} < t \le (1 - q^{-1})n - C_q \sqrt{n \log n}$, 
by~\cref{c:larget} the number of $t$-error correcting codes is at most
\[
\binom{q^n}{\le A_q(n, 2t+1)} \le q^{n A_q(n, 2t+1)} \le 2^{(\log_2 q) \cdot n \cdot o( H_q(n, t)/n)} = 2^{o(H_q(n,t)}. \qedhere
\] 
% \begin{enumerate}[(a)]
% \item This was observed in~\cref{p:thma}.
% \item When $10\sqrt{n} \le t \le (1 - q^{-1})n - C_q \sqrt{n \log n}$, we observe by~\cref{c:larget} and~\cref{eq:weak-ub} that the number of $t$-error correcting codes is at most 
% \[
% q^{n A_q(n, 2t+1)} \le 2^{(\log_2 q) \cdot n \cdot o( H_q(n, t)/n)} = 2^{o(H_q(n,t)}.
% \]
% \item \cref{t:main}(c) follows immediately by~\cref{lem:tc}. \qedhere
%\end{enumerate}
\end{proof}

% \bibliographystyle{amsplain0}
% \bibliography{errorcorrecting.bib}

\begin{thebibliography}{10}

\bibitem{BLS17}
J\'{o}zsef Balogh, Hong Liu, and Maryam Sharifzadeh, \emph{The number of
  subsets of integers with no {$k$}-term arithmetic progression}, Int. Math.
  Res. Not. IMRN (2017), 6168--6186.

\bibitem{BTW16}
J\'{o}zsef Balogh, Andrew Treglown, and Adam~Zsolt Wagner, \emph{Applications
  of graph containers in the {B}oolean lattice}, Random Structures Algorithms
  \textbf{49} (2016), 845--872.

\bibitem{Del73}
P.~Delsarte, \emph{An algebraic approach to the association schemes of coding
  theory}, Philips Res. Rep. Suppl. (1973), vi+97.

\bibitem{FMS20}
Asaf Ferber, Gweneth McKinley, and Wojciech Samotij, \emph{Supersaturated
  sparse graphs and hypergraphs}, Int. Math. Res. Not. IMRN (2020), 378--402.

\bibitem{KW80}
D.~J. Kleitman and K.~J. Winston, \emph{The asymptotic number of lattices},
  Ann. Discrete Math. \textbf{6} (1980), 243--249.

\bibitem{KW82}
Daniel~J. Kleitman and Kenneth~J. Winston, \emph{On the number of graphs
  without {$4$}-cycles}, Discrete Math. \textbf{41} (1982), 167--172.

\bibitem{MS77}
F.~J. MacWilliams and N.~J.~A. Sloane, \emph{The theory of error-correcting
  codes}, North-Holland, 1977.

\bibitem{MRRW77}
Robert~J. McEliece, Eugene~R. Rodemich, Howard Rumsey, Jr., and Lloyd~R. Welch,
  \emph{New upper bounds on the rate of a code via the
  {D}elsarte-{M}ac{W}illiams inequalities}, IEEE Trans. Inform. Theory
  \textbf{IT-23} (1977), 157--166.

\bibitem{Sam15}
Wojciech Samotij, \emph{Counting independent sets in graphs}, European J.
  Combin. \textbf{48} (2015), 5--18.

\bibitem{SAP05}
Alexander Sapozhenko, \emph{Systems of containers and enumeration problems},
  International Symposium on Stochastic Algorithms, Springer, 2005, pp.~1--13.

\bibitem{vL75}
J.~H. van Lint, \emph{A survey of perfect codes}, Rocky Mountain J. Math.
  \textbf{5} (1975), 199--224.

\bibitem{vL99}
J.~H. van Lint, \emph{Introduction to coding theory}, third ed.,
  Springer-Verlag, 1999.

\end{thebibliography}

\end{document}